\documentclass[12pt,twoside]{amsart}
\usepackage{amssymb}

\nonstopmode

\numberwithin{equation}{section}
\font\tengothic=eufm10 scaled\magstep 1
\font\sevengothic=eufm7 scaled\magstep 1
\newfam\gothicfam
      \textfont\gothicfam=\tengothic
      \scriptfont\gothicfam=\sevengothic

\newtheorem{theorem}{Theorem}[section]
\newtheorem{lemma}[theorem]{Lemma}

\newtheorem{corollary}[theorem]{Corollary}

\newtheorem{question}[theorem]{Question}

\theoremstyle{definition}
\newtheorem{definition}[theorem]{Definition} % \theoremstyle{remark}
\newtheorem{remark}[theorem]{Remark}
\newtheorem{example}[theorem]{Example}

\newcommand\codim{\operatorname{codim}}

\newcommand\reg{\operatorname{reg}}

\newcommand{\proj}[1]
{ \mathchoice
           { {\mathbb P}^{#1} }
           { {\mathbb P}^{#1} }
           { {\mathbb P}^{#1} }
           { {\mathbb P}^{#1} }
         }

\newcommand{\acm}{arithmetically Cohen-Macaulay }

\begin{document}
\title[Monomial ideals and the Gorenstein liaison class of a complete
intersection]{Monomial ideals and the Gorenstein liaison class of a complete
intersection}

\author[J.\ Migliore, U.\ Nagel]{J.\ Migliore$^*$, U.\ Nagel$^{**}$} 
%\author[]{}
\address{Department of Mathematics,
        University of Notre Dame, 
        Notre Dame, IN 46556, 
        USA}
\email{Juan.C.Migliore.1@nd.edu}
%\author[]{U.\ Nagel}
\address{Fachbereich Mathematik und Informatik, Universit\"at-Gesamthochschule
Paderborn, D--33095 Paderborn, Germany}
\email{uwen@uni-paderborn.de}

\date{\today}           
\thanks{$^*$ Partially supported  by the Department of Mathematics of the
  University of Paderborn\\
$^{**}$ Partially supported  by the Department of Mathematics of the
  University of Notre Dame }

%%%%%%%%%%%%%%%%%%%%%%%%%%%%%%%%

\begin{abstract}In an earlier work the authors described a mechanism for
lifting monomial ideals to reduced unions of linear varieties.  When the
monomial ideal is Cohen-Macaulay (including Artinian), the corresponding union
of linear varieties is arithmetically Cohen-Macaulay.  The first main result
of this paper is that if the monomial ideal is Artinian then the corresponding
union is in the Gorenstein linkage class of a complete intersection
(glicci).  This technique has some interesting consequences.  For instance,
given any $(d+1)$-times differentiable O-sequence $\underline{H}$, there is a
non-degenerate  \acm reduced union of linear varieties with Hilbert function
$\underline{H}$ which is glicci.  In other words, any Hilbert function that
occurs for \acm schemes in fact occurs among the glicci schemes.  This is not
true for licci schemes.  Modifying our technique, the second main result is
that any Cohen-Macaulay  Borel-fixed monomial ideal is glicci.  As a
consequence, all \acm subschemes of projective space are glicci up to flat
deformation.  
\end{abstract}

%%%%%%%%%%%%%%%%%%%%%%%%%%%%%%%

\maketitle

\tableofcontents

%%%%%%%%%%%%%%%%%%%%%%%%%%%%%%%%%%%%%%%%%%%%%%%%%

 \section{Introduction} \label{intro} 

Liaison theory has reached a very satisfying state in codimension two, but in
higher codimension there are still many open problems.  Much of the theory has
been built around linking with complete intersections, called CI-liaison
theory.  However, it has long been known that more generally it is also
possible to link using arithmetically Gorenstein schemes (cf.\
\cite{schenzel} for example).  Indeed, a development of G-liaison theory is
possible (cf.\ \cite{KMMNP}, \cite{migliore}, \cite{nagel}).  In practice,
however, this has been studied less because it is not so easy to find
arithmetically Gorenstein schemes other than complete intersections,
especially containing a given scheme.   Note that in codimension two all
arithmetically Gorenstein schemes are complete intersections, so both
theories include the codimension two case.  

Nevertheless, there has been recent work in the direction of G-liaison theory,
most notably in \cite{KMMNP} where a very geometric approach is taken and
where this theory is compared and contrasted with the more classical
CI-liaison theory.  See also \cite{migliore} for extensive background and
comparisons.  

In the codimension two case, one of the first important results was the
theorem of Gaeta that every arithmetically Cohen-Macaulay, codimension two
scheme is in the liaison class of a complete intersection (i.e.\ is {\em
licci}, a term introduced in \cite{HU}).  In \cite{KMMNP} the authors
introduced the notion of {\em glicci} schemes, i.e.\ those which are in the
Gorenstein liaison class of a complete intersection.  They generalized
Gaeta's theorem by showing that every scheme which arises as the maximal
minors of a homogeneous matrix (and which have the right codimension
depending on the size of the matrix) is glicci.  (Note that \acm schemes of
codimension two satisfy this property, thanks to the Hilbert-Burch theorem.) 
However, the authors of \cite{KMMNP} asked if a more general result might
hold:

\begin{question}[\cite{KMMNP}] \label{mainq} Is it true that all \acm
subschemes of $\proj{n}$ are glicci?
\end{question}

  Some evidence of this was provided by showing that on a smooth
rational \acm surface in $\proj{4}$ all \acm curves (i.e.\ divisors) are
glicci.  Casanellas and Mir\'o-Roig \cite{CM1} extended this by finding a
large class of smooth surfaces in $\proj{4}$ where the same conclusion holds,
and in a more recent paper \cite{CM2} they extend this to a large class of
smooth schemes of any codimension.

  In this paper we make some further progress in this direction.  We prove
the glicciness of two different kinds of Cohen-Macaulay ideals.  First we
recall that if $J$ is any Artinian monomial ideal then it is shown in
\cite{MN2} how to produce the ideal $I$ of a non-degenerate \acm reduced union
of linear varieties of any dimension whose Artinian reduction is precisely
$J$.  The first main result of this paper is that any such $I$ is glicci.  
As a corollary we get that given any numerical function which occurs as the
Hilbert function of some non-degenerate \acm subscheme of $\proj{n}$ of any
codimension, there is a reduced, glicci subscheme with precisely that Hilbert
function.  Example \ref{false-for-licci} shows that this is not true if we
replace Gorenstein links by complete intersection links.

Our second main result is that any Cohen-Macaulay Borel-fixed monomial ideal
(Artinian or not) is glicci.  This result is of a rather general nature. 
Indeed, it is well-known that every generic initial ideal of an \acm
subscheme is a Cohen-Macaulay Borel-fixed ideal which defines a deformation
of the original scheme.  Thus our result says that every \acm subscheme
admits a flat deformation which is glicci.  In other words, we have found an
affirmative answer to Question \ref{mainq} ``up to flat deformation.''

%%%%%%%%%%%%%%%%%%%%%%%%%%%%%%%%%%%%%%%%%%%%%%%%%%%%%%%%%%%%%%%%%%%%%%%%%%%%

\section{Preliminaries}

Let $K$ be an infinite field and let $S = K[x_1,\dots,x_n]$ and
$R=K[x_1,\dots,x_n,u_1,\dots,u_t]$, $t \geq 0$.  We first recall the set-up
and one of the main results of \cite{MN2}.  

\begin{definition}\label{t-lifting-def}
Let $I \subset R$ and $J \subset S$ be homogeneous ideals.  Then we say $I$ is
a {\em $t$-lifting} of $J$ to $R$ (or when $R$ is understood, simply a {\em
$t$-lifting} of $J$) if $(u_1,\dots,u_t)$ is a regular sequence
on $R/I$ and $(I,u_1,\dots,u_t)/(u_1,\dots,u_t) \cong J$.  
\end{definition}

The definition of a $t$-lifting can be extended to modules, but Definition
\ref{t-lifting-def} suffices for our purposes.  Consider now a matrix of
linear forms
\[
A = 
\left [
\begin{array}{ccccccccc}
L_{1,1} & L_{1,2} & L_{1,3} & \dots \\
L_{2,1} & L_{2,2} & L_{2,3} & \dots \\
\vdots & \vdots & \vdots \\
L_{n,1} & L_{n,2} & L_{n,3} & \dots 
\end{array}
\right ]
\]
where the $L_{j,i}$ are in $R$.  $A$ will be called the {\em lifting matrix},
for reasons that will be apparent shortly.  For now we assume that there are
infinitely many columns, but in practice when we have a specific ideal $J 
\subset S$ that we want to lift we can assume that the number of columns is
finite, for instance equal to the regularity of $J$ (or less).  Assume that
the polynomials $F_j = \prod_{i=1}^N L_{j,i}$, $1 \leq j \leq n$, define a 
complete intersection, $X$.  Note that $F_j$ is the product of the entries of
the $j$-th row, and that the height of the complete intersection is $n$, the
number of variables in $S$.

Let $m = \prod_{j=1}^n x_j^{a_j} \subset S$ be a monomial.  We associate to
$m$ the homogeneous polynomial
\[
\bar m = \prod_{j=1}^n \left ( \prod_{i=1}^{a_j} L_{j,i} \right ) \in R.
\]
Let $J = (m_1,\dots,m_r) \subset S$ be a monomial ideal.  Associated to $J$
we define the ideal $I = (\bar m_1,\dots, \bar m_r) \subset R$.

\begin{theorem}[\cite{MN2}] \label{MN2-main}
\begin{itemize}
\item[(i)] The ideal $I$ is saturated.

\item[(ii)] $S/J$ is Cohen-Macaulay (including the case where it is Artinian)
if and only if $R/I$ is Cohen-Macaulay.  In fact, $I$ (as an $R$-module) and
$J$ (as an $S$-module) have the same graded Betti numbers.

\item[(iii)] If, for each $j$ and $i$, we have $L_{j,i} \in
K[x_j,u_1,\dots,u_t]$ then $I$ is a $t$-lifting of $J$.  Otherwise we say
that $I$ is a {\em pseudo-lifting} of $J$.
\end{itemize}
\end{theorem}

If the entries of $A$ are chosen sufficiently generally then $I$ in fact
defines a reduced union of linear varieties with good intersection
properties.  Now we recall some results from \cite{KMMNP}.

We now recall the notion of
{\em Basic Double G-linkage} introduced in \cite{KMMNP}, so called because of
part (iv) and the notion of Basic Double Linkage (\cite{lazarsfeld-rao},
\cite{BM4}, \cite{GM4}).

\begin{theorem}[\cite{KMMNP} Lemma 4.8, Remark 4.10 and Proposition
5.10]\label{KMMNP-lemma} Let $J \supset I$ be homogeneous ideals of $R' =
K[x_0,\dots,x_n]$, defining schemes $W \subset V \subset \proj{n}$ such that
$\codim W = \codim V + 1$.  We also allow the possibility that $J$ is
Artinian and $V$ is a zeroscheme.  Let $A \in R'$ be an element of degree
$d$ such that $I:A = I$.  Then we have
\begin{itemize}
\item[(i)] $\deg (I+A \cdot J) = d \cdot \deg I + \deg J$.

\item[(ii)] If $I$ is perfect and $J$ is unmixed then $I+A\cdot J$ is unmixed.

\item[(iii)] $J/I \cong [(I+A\cdot J)/I] (d).$

\item[(iv)] If $V$ is \acm and generically Gorenstein and $J$ is unmixed then
$J$ and $I + A\cdot J$ are linked in two steps using Gorenstein ideals.

\item[(v)] The Hilbert functions are related by 
\[
\begin{array}{rcl}
h_{R'/(I+A\cdot J)} (t) & = & h_{R'/(I+(A))} (t) + h_{R'/J}(t-d) \\
& = &h_{R'/I}(t) - h_{R'/I}(t-d) + h_{R'/J}(t-d)
\end{array}
\]

\end{itemize}
\end{theorem}

Theorem \ref{KMMNP-lemma} should be interpreted as viewing the scheme $W$
defined by $J$ as a divisor on the scheme $V$ defined by $I$, and adding to it
a hypersurface section $H_A$ of $V$ defined by the polynomial $A$.  Note
that $I_{H_A} = I_V + (A)$.  If $V$ and $W$ are \acm then the divisor $W+H_A$
is again \acm (by step 4).  As an immediate application we have the following
by successively applying Theorem \ref{KMMNP-lemma}.

\begin{corollary}[\cite{MN3}] \label{hyperplane-sects}
Let $V_1 \subset V_2 \subset \dots \subset V_r \subset \proj{n}$ be \acm
schemes of the same dimension, all generically  Gorenstein. 
Let
$H_1,\dots,H_r$ be hypersurfaces, defined by forms $F_1,\dots,F_r$, such that
for each $i$,
$H_i$ contains no component  of $V_j$ for any $j \leq i$.  Let $W_i$ be the
\acm schemes defined by the corresponding hypersurface sections: $I_{W_i} =
I_{V_i} + (F_i)$.  Then we have the following.

\begin{itemize}
\item[(i)] Viewed as divisors on $V_r$, the sum $Z$ of the $W_i$ (which
is just the union if the hypersurfaces are general enough) is in the same
Gorenstein liaison class as $W_1$.

\item[(ii)] In particular, $Z$ is arithmetically Cohen-Macaulay.  

\item[(iii)] As ideals we
have
\[
I_Z = I_{V_r} + F_r \cdot I_{V_{r-1}} + F_r F_{r-1} I_{V_{r-2}} + \cdots
+ F_r F_{r-1} \cdots F_2 I_{V_1} + (F_r F_{r-1}\cdots F_1).
\]

\item[(iv)] Let $d_i = \deg F_i$.  The Hilbert functions are related by the
formula
\[
\begin{array}{rcl}
h_Z (t) & = & h_{W_r} (t) + h_{W_{r-1}} (t-d_{r}) + h_{W_{r-2}}
(t-d_{r}-d_{r-1}) +
\dots \\
&& + h_{W_1} (t-d_{r}-d_{r-1}-\dots-d_2).
\end{array}
\]
\end{itemize}
\end{corollary}

\begin{corollary}\label{glicci-cor} We keep the notation of Corollary
\ref{hyperplane-sects}.  If $V_1$ is glicci then so is $Z$.  
\end{corollary}

\begin{proof}
This follows from part (iv) of Theorem \ref{KMMNP-lemma}, and from the fact
that Gorenstein liaison is preserved under hypersurface sections
(\cite{migliore} Proposition 5.2.17).  Note that the reverse direction is not
necessarily true (or in any case is not known to be true): if $W_1$ is
glicci, it does not necessarily (or at least immediately) hold that $V_1$
is.  See 
\cite{migliore} Example 5.2.26 for some discussion.
\end{proof}

We now discuss the decomposition of a monomial ideal, which we will use in
the remaining sections. 

\begin{definition} \label{def-of-LS-BF}
Let $>$ denote the degree-lexicographic order on monomial ideals, i.e.\
$x_1^{a_1}\cdots x_n^{a_n} > x_1^{b_1}\cdots x_n^{b_n}$ if the first nonzero
coordinate of the vector
\[
\left ( \sum_{i=1}^n (a_i - b_i), a_1 - b_1 ,\dots,a_n - b_n \right )
\]
is positive.  Let $J$ be a monomial ideal.  Let $m_1,m_2$ be monomials in
$S$ of the same degree such that $m_1 > m_2$.  Then $J$ is a {\em lex-segment
ideal} if $m_2 \in J$ implies $m_1 \in J$.  When $\hbox{char}(K) = 0$, we say
that
$J$ is a {\em Borel-fixed ideal} if 
\[
m = x_1^{a_1}\cdots x_n^{a_n} \in J, \ a_i >0, \hbox{ implies }
\frac{x_j}{x_i} \cdot m \in J
\]
for all $1 \leq j < i \leq n$.
\end{definition}

\begin{remark}\label{BF-rmk}
Definition \ref{def-of-LS-BF} says that if $J$ is Borel-fixed and $m \in J$
is a monomial then one can reduce any power of a variable occurring in $m$ by
one and increase the power of a larger variable by one, and the result is
again in $J$.  Note that this is not the same as lex-segment.  For example,
in the ring $K[x_1,x_2,x_3 ]$ consider the ideal
\[
J = \langle x_1^3, x_1^2x_2, x_1x_2^2 \rangle.
\]
This is Borel-fixed but not lex-segment, since $x_1^2x_3 \notin J$.  The two
notions are not even equivalent in the Artinian case, as the same example
shows if we adjoin to $J$ all monomials of degree 4.  However, a lex-segment
ideal is always Borel-fixed.
\end{remark}

\begin{lemma}\label{decomp}
Let $J \subset S=K[x_1,\dots,x_n]$ be a monomial ideal.  Let $\alpha$ be the
highest power of $x_1$ occurring in a minimal generator of $J$.  Then there
is a  uniquely determined decomposition
\[
J = \sum_{j=0}^\alpha x_1^j \cdot (I_j \cdot S)
\]
where $I_0 \subset I_1 \subset \cdots \subset I_{\alpha-1} \subset
I_{\alpha} $ are monomial ideals in $T = K[x_2,\dots,x_n]$.   Furthermore, 
\begin{itemize}
\item[(i)] $I_j = (J : x_1^j) \cap K[x_2,\dots,x_n]$. 
\item[(ii)] If $J$ is Artinian then so is each $I_j$ and $I_\alpha = (1)$.

\item[(iii)] Assume $\hbox{char}(K) = 0$.  If $J$ is a Borel-fixed ideal
(e.g.\ a lex-segment ideal) then
 $\alpha$ is the initial degree of $J$, $I_\alpha = (1)$, and each $I_j$ is
again Borel-fixed.
\end{itemize}
\end{lemma}

\begin{proof}
The case of Artinian lex-segment ideals was observed in \cite{MN3}. 

The existence of the decomposition is clear if we choose the ideals $I_j$ as 
described in (i). Conversely, if we have the decomposition, then we get in
case  $0 \leq j \leq \alpha$: 
\[
\begin{array}{rcl} 
J : x_1^j & = & \displaystyle  \sum_{k=0}^j I_k \cdot S + 
 \sum_{k=j+1}^\alpha x_1^{k-j}  \cdot (I_k \cdot S) \\
& = & \displaystyle I_j \cdot S +   \sum_{k=j+1}^\alpha x_1^{k-j} \cdot (I_k
\cdot S), 
\end{array} 
\] 
thus 
$$
(J : x_1^j) \cap T = (I_j \cdot S) \cap T = I_j
$$ 
proving (i) and the uniqueness of the decomposition. 

For (ii), since $J$ is Artinian then it contains pure
powers of $x_2,\dots,x_n$, so these are automatically in $I_0$, making $I_0$
Artinian.  Then the inclusions imply that the other $I_j$ are also Artinian. 
Furthermore, $x_1^\alpha$ is a minimal generator of $J$, so $I_\alpha = (1)$
as claimed.  For part (iii), the hypothesis implies that $x_1^\alpha$ is a
minimal generator of $J$.  The fact that $I_j$ is Borel-fixed follows
immediately from the definition of Borel-fixed and the description of $I_j$
in the statement of the Lemma.  
\end{proof}

\begin{lemma} \label{hf-lemma}
Keeping the notation of Lemma \ref{decomp}, for any $s\geq 0$, we have
\[
h_{S/J}(s) = \sum_{j=0}^{\alpha-1} h_{T/I_j}(s-j) + h_{S/I_{\alpha} \cdot
S}(s-  \alpha). 
\]
\end{lemma}

\begin{proof}
If $\alpha = 0$ then $J = I_0 \cdot S$ and the claim is clear. If $\alpha > 0$
then  multiplication by $x_1$ provides the exact sequence 
$$
0 \to S/\sum_{j=1}^{\alpha} x_1^{j-1} (I_j \cdot S) (-1) \stackrel{\cdot
x_1}{\longrightarrow} S/J \to T/I_0 \to 0. 
$$
Hence the claim follows by induction on $\alpha$. 
\end{proof} 

\begin{remark}
If $J$ is Artinian or Borel-fixed then the Hilbert function formula of Lemma
\ref{hf-lemma} simplifies to
\[
h_{S/J}(s) = \sum_{j=0}^{\alpha-1} h_{T/I_j}(s-j)
\]
since in either of these cases $I_\alpha = (1)$ by Lemma \ref{decomp}.
\end{remark}

%%%%%%%%%%%%%%%%%%%%%%%%%%%% section 3 %%%%%%%%%%%%%%%%%%%%%%%%%%%%%%%%%%%

\section{Glicci Ideals}

Let $J$ be an Artinian monomial ideal in $S = K[x_1,\dots,x_n]$. 
Let $A$ be a lifting matrix for $J$ and
assume that the entries of $A$ are sufficiently general so that the
lifted ideal is a reduced union of linear varieties..  The number of columns 
of $A$ only has to be as large as the largest degree of a minimal generator
of $J$; if $J$ is lex-segment then this degree is 
$= \reg (J)$.  Applying the pseudo-lifting procedure
described in Section 2, we get an ideal $I \subset R =
K[x_1,\dots,x_n,u_1,\dots,u_t]$ which, by Theorem
\ref{MN2-main}, is the saturated ideal of an \acm subscheme $Z$ of
$\proj{n+t-1}$ of codimension $n$.

\begin{theorem} \label{main-th}
$Z$ is glicci.
\end{theorem}

\begin{proof}
The proof is by induction on $n$, the codimension.  For codimension two it is
known that any \acm subscheme of projective space is licci, so there is
nothing to prove.  Hence we assume $n \geq 3$.

By Lemma \ref{decomp} we have
\begin{equation}\label{Jdecomp}
J = \sum_{j=0}^{\alpha} x_1^j \cdot I_j
\end{equation}
where $I_0 \subset I_1 \subset \cdots \subset I_{\alpha-1} \subsetneq
I_{\alpha} = S$ and for each $j$, $I_j$ is an Artinian  ideal in
$K[x_2,\dots,x_n]$.  Notice that the lifting matrix $A$ has $n$ rows, and if
we remove the first row then the remaining matrix $A'$ can be used to lift the
ideals $I_j$.

Let $\bar I_j$ be the ideal obtained by lifting $I_j$ using $A'$.  Let $Y_j$
be the \acm subscheme of $\proj{n+t-1}$ defined by $\bar I_j$.  Note that
$Y_j$ has codimension $n-1$, but the projective space does not change since
the linear forms which are the entries of $A$ were taken from the ring $R$. 
Note also that $Y_{\alpha-1} \subset \cdots \subset Y_1 \subset Y_0$ are \acm
schemes of the same dimension and generically Gorenstein, as in the set-up of
Corollary \ref{hyperplane-sects}.

Thanks to (\ref{Jdecomp}) we have
\[
I = \bar I_0 + L_{1,1} \cdot \bar I_1 + L_{1,1} L_{1,2} \cdot \bar I_2 +
\cdots + L_{1,1} L_{1,2} \cdots L_{1,\alpha-1} \cdot \bar I_{\alpha-1} +
(L_{1,1}L_{1,2} \cdots L_{1,\alpha-1} L_{1,\alpha} )
\]
Hence by Corollary \ref{hyperplane-sects} (iii), the scheme $Z$ obtained from
lifting is in fact also obtained by taking the union of the successive
hypersurface sections of the $Y_j$.  By induction, $Y_{\alpha-1}$ is glicci.
By Corollary \ref{glicci-cor}, then, $Z$ is also glicci.
\end{proof}

As a corollary of Theorem \ref{main-th} we would like to show that given
``any'' Hilbert function we can find a glicci subscheme with that Hilbert
function.  Recall from \cite{GMR} that the Hilbert functions which can occur
for \acm subschemes of a given dimension $d$ have been completely
characterized.  Indeed, for a function $f : {\mathbb Z} \rightarrow {\mathbb
Z}$ we define the first difference $\Delta f$ by $\Delta f (n) = f(n) -
f(n-1)$ and the $k$-th difference $\Delta^k f$ by iteration.  An O-sequence
is one that satisfies Macaulay's growth condition \cite{mr.macaulay}.  A
$k$-times differentiable O-sequence is one for which also the first $k$
differences are O-sequences.  Then a function $f:{\mathbb N} \rightarrow
{\mathbb N}$ occurs as the Hilbert function of some $d$-dimensional \acm
scheme (in fact it can always be chosen reduced) if and only if $f$ is a
$(d+1)$-times differentiable O-sequence.

We immediately get the following somewhat surprising conclusion.

\begin{corollary} \label{hilbftn}
Let $\underline{H}$ be any $(d+1)$-times differentiable
O-sequence.  Then $\underline{H}$ occurs as the Hilbert function of some
non-degenerate glicci subscheme of projective space.
\end{corollary}

\begin{proof}
Let $\underline{h}$ be the $(d+1)$-st difference of $\underline{H}$.  Let $J$
be the Artinian lex-segment ideal with Hilbert function $\underline{h}$.  If
$J \subset K[x_1,\dots,x_n]$ then choose a lifting matrix with entries
$L_{j,i}
\in K[x_j,u_1,\dots,u_{d+1}]$.  The lifted ideal $I$ defines a glicci
subscheme of $\proj{n+d-1}$ by Theorem \ref{main-th}, and it has Hilbert
function $\underline{H}$ since it is a $(d+1)$-lifting.  The non-degenerate
property comes directly from the lifting, cf.\ \cite{MN2}.
\end{proof}

\begin{example}\label{false-for-licci}
We remark that Corollary \ref{hilbftn} is false for complete intersection
liaison.  Indeed, the $h$-vector $(1,3)$ cannot occur for any codimension 3
licci subscheme of projective space.  To see this, note that the minimal free
resolution of any \acm subscheme with this $h$-vector is linear, and 
\cite{HU}, Corollary 5.13, then guarantees that it is not licci.  (Note that
degenerate subschemes of projective space, of codimension $>3$, could also
have this $h$-vector, and we do not know if the ``extra room'' makes a
difference in the non-licciness.)
\end{example}

From the proof of Corollary \ref{hilbftn} one would be very tempted to
conclude that we have proved that any Artinian monomial ideal is glicci,
since liaison is preserved under general hyperplane sections, even for the
Artinian reduction (cf.\ \cite{migliore} Remark 5.2.18).  However, the
proofs above use {\em bilinks}, so even if a variable $u_i$ is a non
zero-divisor for the scheme $Z$ (and hence any of its components), it is not
necessarily true that the same true for the linked schemes.  However, we can
obtain an important case of this result, and in fact more, by modifying the
above approach slightly.

From now on we assume $\hbox{char} (K) = 0$.  We begin with a lemma.

\begin{lemma} \label{BF-equiv}
Let $J$ be a Borel-fixed monomial ideal of codimension $c$.  The following
are equivalent.
\begin{itemize}
\item[(i)] $J$ is Cohen-Macaulay.
\item[(ii)] $J$ is equidimensional.
\item[(iii)] $J$ contains a pure power of $x_c$, and the
variables $x_{c+1},\dots,x_n$ do not occur in any of the minimal generators.
\item[(iv)] $J$ is a cone over an Artinian Borel-fixed ideal in
$K[x_1,\dots,x_c]$.
\end{itemize}
\end{lemma}

\begin{proof}
The implications (i) $\Rightarrow$ (ii) and (iv)
$\Rightarrow$ (i) are always true. Note that condition (iii) implies that 
$J$ in fact contains pure powers of each of the variables $x_1,\dots,x_c$, by
the Borel-fixed property.  Then the implication (iii) $\Rightarrow$ (iv)
is immediate, since Borel-fixed is already assumed.  

 So we have only to prove (ii) $\Rightarrow$ (iii).  Since $J$ has
codimension $c$, it contains a regular sequence of length $c$.  By the
Borel-fixed property we may take this regular sequence to consist of pure
powers of variables, and again by the Borel-fixed property we can take it to
be powers of $x_1,\dots,x_c$.  Suppose that one of the other variables, say
$x_{c+1}$, occurs in one of the minimal generators of $J$ to some power $a
\geq 1$.  By a standard trick on monomial ideals (cf.\ for instance
\cite{eisenbud} Exercise 3.8) we can then decompose $J$ as $J = A \cap (J +
(x_{c+1}^a))$ where $A$ is again a monomial ideal.  But this shows that the
primary decomposition of $J$ has at least one component of height $c+1$,
contradicting the hypothesis that $J$ is equidimensional.
\end{proof}

\begin{theorem}\label{BF-glicci}
Any Cohen-Macaulay Borel-fixed monomial ideal is glicci.
\end{theorem}

\begin{proof}
Let $J$ be a Cohen-Macaulay Borel-fixed monomial ideal in $S
= K[x_1,\dots,x_n]$ of height $c$.  By Lemma \ref{BF-equiv}, we may view $J$
as a cone over an Artinian Borel-fixed ideal in $K[x_1,\dots,x_c]$.  By Lemma
\ref{decomp}, 
\[
J = \sum_{j=0}^\alpha x_1^j \cdot I_j
\]
where $\alpha$ is the initial degree of $J$ and the $I_j$ are cones over 
Artinian Borel-fixed ideals in $K[x_2,\dots,x_c]$ satisfying $I_0 \subset
I_1 \subset \dots$.  We can rewrite this as
\[
J = I_0 + x_1 \cdot I'
\]
where $I_0$ is a cone over an Artinian Borel-fixed ideal in 
$K[x_2,\dots,x_c] \subset S$, and $I'$ is a Borel-fixed monomial ideal in $S$
whose initial degree is one less than that of $J$.  

Following Theorem
\ref{MN2-main}, we can lift $I_0$ to an ideal $\bar I_0$ in
$K[x_1,\dots,x_c] \cap S$; that is, we choose a lifting matrix $A$ whose
entries are linear forms
$L_{j,i} \in K[x_j,x_1]$, $2 \leq j \leq n$.  For example, take
\[
A = 
\left [ 
\begin{array}{ccccccccc}
x_2 & x_2 + x_1 & x_2 + 2 x_1 & x_2 + 3 x_1  \dots \\
x_3 & x_3 + x_1 & x_3 + 2 x_1 & x_3 + 3 x_1  \dots \\
\vdots & \vdots & \vdots & \vdots \\
x_n & x_n + x_1 & x_n + 2 x_1 & x_n + 3 x_1  \dots 
\end{array}
\right ].
\]

We now make some observations.

\begin{itemize}

\item[(1)] $I'$ is also
Cohen-Macaulay by Lemma \ref{BF-equiv}, and it has the same height $c$ as $J$
since it contains a complete intersection consisting of powers of
$x_1,\dots,x_c$.

\item[(2)] $I_0 \subset I'$.  This follows from the sequence of inclusions on
the $I_j$.

\item[(3)] $\bar I_0 \subset I'$.  This follows immediately.  For
instance, suppose that $x_2^3 x_3^4 \in I_0$.  Then
\[
x_2 (x_2 + x_1)(x_2 + 2x_1)(x_3)(x_3+x_1)(x_3+2x_1)(x_3+3x_1) \in \bar I_0.
\]
By the Borel-fixed property of $J$ and the fact that $I_0 \subset I'$, it
follows immediately that each term of this polynomial is in $I'$.

\item[(4)] $I_0$ and $\bar I_0$ are both Cohen-Macaulay, and $\hbox{ht}( \bar
I_0) = \hbox{ht}( I_0) = c-1$.  This  follows from the fact that $I_0$ is
Cohen-Macaulay by Lemma \ref{BF-equiv} and that the Cohen-Macaulay property
and the codimension are preserved under lifting.
\end{itemize}

Let $\bar J = \bar I_0 + x_1 \cdot I'$.   An analysis similar to observation 
(3) above shows quickly that $\bar J \subset J$.  But both are Cohen-Macaulay
of the same height in $S$, and they have the same Hilbert function (since the
Hilbert function of $K[x_2,\dots,x_c]/I_0$ is the first difference of that of
$K[x_1,\dots,x_c]/\bar I_0$).  Hence we obtain that $\bar J = J$.

Although $I_0$ is not necessarily generically Gorenstein, the lifting results
guarantee that $\bar I_0$ is, and we have noted that $\bar I_0$ is
Cohen-Macaulay.  Hence Theorem \ref{KMMNP-lemma} (iv) says that $J = \bar J$
is G-bilinked to $I'$, and in particular $I'$ is Cohen-Macaulay.  We have
noted that the initial degree of
$I'$ is one less than that of $J$.  Hence in a finite (even) number of steps
we obtain that $J = \bar J$ is linked to the hyperplane section $\bar I_0 +
(x_1) = I_0 + (x_1)$.  Thus it is enough to show that $I_0$ is glicci. 
Let $J_0$ denote the ideal $I_0
\cap T$ in $T :=K[x_2,...,x_n]$.  Then $I_0$ is just a cone over $J_0$. 
By induction on the height, $J_0$ is glicci in $T$.  Then taking  cones we
get that also $I_0$ is glicci. Hence we
have shown that $J = \bar J$ is glicci, as claimed.
\end{proof}

\begin{remark}
Theorem \ref{BF-glicci} is of a rather general nature.  It is
well-known that every generic initial ideal of an \acm subscheme is a
Cohen-Macaulay Borel-fixed ideal which defines a deformation of the original
scheme.  Indeed, the fact that it is Borel-fixed is due to Galligo
\cite{galligo};  that it gives a flat deformation is due to Bayer
\cite{bayer}; that it is again Cohen-Macaulay follows from a result of Bayer
and Stillman (cf.\ \cite{eisenbud} Theorem 15.13).  Thus our result says that
every \acm subscheme admits a flat deformation which is glicci.  In other
words, we have found an affirmative answer to Question \ref{mainq} ``up to
flat deformation.''
\end{remark}

\begin{example}
We illustrate the above ideas by finding a glicci subscheme $Z \subset
\proj{3}$ with $h$-vector 
\[
\underline{h} = (1,3,6,10,4,2).
\]
  Note that using complete intersections it does not seem
promising that a licci subscheme with this $h$-vector can be found since the
smallest complete intersection containing it would be the complete
intersection of three quartics, and the residual would have even larger degree
and will not lie in a smaller complete intersection.

Instead we consider the ring $S = K[x_1,x_2,x_3]$ and let $J$ be the Artinian
lex-segment ideal with Hilbert function $\underline{h}$.  We have the
decomposition 
\[
J = I_0 + x_1 \cdot I_1 + x_1^2 \cdot I_2 + x_1^3 \cdot I_3 + (x_1^4),
\]
 where the
$I_j$ are Artinian lex-segment ideals in $T = K[x_2,x_3]$ whose Hilbert
functions are given as follows (note the shifting for $I_1$, $I_2$ and $I_3$
to apply Lemma \ref{hf-lemma}):

\smallskip

\hskip 2in
\begin{tabular}{r|cccccccccc}
&\multicolumn{3}{l}{degree:} \\
& 0 & 1 & 2 & 3 & 4 & 5 & 6 \\ \hline
$I_0$ & 1 & 2 & 3 & 4 & 4 & 2 \\
$I_1$ && 1 & 2 & 3 \\
$I_2$ &&& 1 & 2 \\
$I_3$ &&&& 1 \\ \hline
& 1 & 3 & 6 & 10 & 4 & 2
\end{tabular}

\bigskip

If 
\[
A = 
\left [
\begin{array}{ccccccccc}
L_{1,1} & L_{1,2} & L_{1,3} & \dots \\
L_{2,1} & L_{2,2} & L_{2,3} & \dots \\
L_{3,1} & L_{3,2} & L_{3,3} & \dots 
\end{array}
\right ]
\]
 is a lifting matrix with 3 rows and at least 6 columns then the lifted
ideal $I$ is the saturated ideal of a zeroscheme $Z$ in $\proj{3}$ which 
\begin{itemize}
\item[(i)] is reduced if $A$ is sufficiently general, 

\item[(ii)] is glicci, by Theorem \ref{main-th}, and  

\item[(iii)] has $h$-vector $\underline{h}$.  
\end{itemize}

The proof of Theorem
\ref{main-th} shows that $Z$ can in fact be obtained as the union of
successive hyperplane sections (denoting hyperplanes with the same notation
as the corresponding linear forms)
\[
Z = (L_{1,1}\cap V_0) \cup (L_{1,2} \cap V_1) \cup (L_{1,3} \cap V_2)
\cup (L_{1,4} \cap V_3)
\]
where
\[
V_3 \subset V_2 \subset V_1 \subset V_0
\]
are reduced \acm configurations of lines in $\proj{3}$ obtained by lifting
$I_0,\dots,I_3$ using the submatrix
\[
A' = 
\left [
\begin{array}{ccccccccc}
L_{2,1} & L_{2,2} & L_{2,3} & \dots \\
L_{3,1} & L_{3,2} & L_{3,3} & \dots 
\end{array}
\right ]
\]
and the $h$-vectors of the $V_j$ are given by the rows of the table above.
\end{example}

%%%%%%%%%%%%%%%%%%%%%%%%%%%%%%%%%%%%%%%%%%%%%%%%%%%%%%%%%%%%%%%%%%%%%%%%%

\section{Further Comments} 

We end with some comments and questions raised by this paper.  The results in
this paper, as well as those in \cite{KMMNP}, \cite{CM1} and \cite{CM2}, 
suggest strongly to us that the answer to Question \ref{mainq} is ``yes.'' 
The following ideas may help to ultimately give a final answer to this
question.

\begin{enumerate}
\item We have seen that Cohen-Macaulay Borel-fixed monomial ideals are
glicci.  Is it in fact true that {\em every} Cohen-Macaulay monomial ideal is
glicci?  Or is it at least true that every Artinian monomial ideal is
glicci?

\item Given a Hilbert function, our lifting gives
the ``worst'' \acm scheme with that Hilbert function.  As a result, this
scheme should be the most difficult to find ``good'' arithmetically
Gorenstein  schemes containing it.  Since we can find suitable ones for this
``worst case,'' can this suggest how to link a different \acm scheme with
that same Hilbert function down to a complete intersection?
\end{enumerate}

\end{document}